\documentclass[12pt]{article}
\usepackage{amsmath, amsthm, amssymb}
\numberwithin{equation}{section}
\begin{document}
\author{Ajai Choudhry and Oliver Couto}
\title{A new diophantine equation\\ involving fifth powers}
\date{}
\maketitle

\begin{abstract}
In this paper we obtain a parametric solution of the hitherto unsolved diophantine equation $(x_1^5+x_2^5)(x_3^5+x_4^5)=(y_1^5+y_2^5)(y_3^5+y_4^5)$. Further, we show, using elliptic curves, that there exist infinitely many parametric solutions of the aforementioned diophantine equation, and they can be effectively computed.
\end{abstract}



\section{Introduction}
In this paper we obtain a parametric solution of the following diophantine equation involving fifth powers:
\begin{equation}
(x_1^5+x_2^5)(x_3^5+x_4^5)=(y_1^5+y_2^5)(y_3^5+y_4^5). \label{eq5}
\end{equation}

Eq. \eqref{eq5} has not been considered at all in the existing literature. We note that when $n$ is a positive integer $ < 5$, parametric solutions of the diophantine equation,
\begin{equation}
(x_1^n+x_2^n)(x_3^n+x_4^n)=(y_1^n+y_2^n)(y_3^n+y_4^n). \label{eqn}
\end{equation}
are readily obtained since we know complete/ parametric solutions of the diophantine equation $x_1^n+x_2^n=y_1^n+y_2^n$. In fact, it is not difficult to obtain the complete solution in rational numbers of Eq. \eqref{eqn} when $n=2$ or $n=3$. However, no nontrivial integer solutions of the equation $x_1^5+x_2^5=y_1^5+y_2^5$ are known, and there is no obvious way of obtaining nontrivial integer solutions of Eq. \eqref{eq5}.

We also note that parametric solutions of certain  diophantine equations involving fifth powers are already known. For instance, parametric solutions of the diophantine equations,
\begin{equation}
\sum_{i=1}^sx_i^5=\sum_{i=1}^sy_i^5, \quad {\rm where} \;\;s=3 \;\;{\rm or}\;\; s=4, 
\end{equation}
\begin{equation}
ax_1^5+bx_2^5+cx_3^5=ay_1^5+by_2^5+cy_3^5, 
\end{equation}
 where $a, b, c$ are distinct nonzero integers such that $a+b+c = 0$, and 
\begin{equation}
ax_1^5+bx_2^5+cx_3^5+dx_4^5=ay_1^5+by_2^5+cy_3^5+dy_4^5, 
\end{equation}
 where $a, b, c, d$ are arbitrary nonzero integers, have been given by Bremner \cite{Br},  Moessner \cite{Mo}, Swinnerton-Dyer \cite{Sw},  Lander \cite{La} and  by  Choudhry \cite{Ch1, Ch2, Ch3}.  In all these cases, the degree of the diophantine equation is 5. In contrast, however,  Eq. \eqref{eq5} is of degree 10, and till now there are no instances of integer solutions of such a diophantine equation being published.

\section{The diophantine equation \eqref{eq5}}\label{diopheqneq5}

\subsection{Preliminary remarks}\label{prelim} 
There are certain integer solutions of Eq. \eqref{eq5} that satisfy Eq. \eqref{eqn}  when $n$ is any arbitrary odd positive  integer. All such solutions of Eq. \eqref{eq5} will be considered as trivial solutions. An example of such a trivial solution is 
\begin{equation}
\begin{aligned}
x_1 &=& a_1u, \quad & x_2 &=& a_2u, \quad & x_3 &=& a_3v, \quad & x_4 &=& a_4v,\\
 y_1 &=& a_1v, \quad & y_2 &=& a_2v, \quad & y_3 &=& a_3u, \quad & y_4 &=& a_4u,
\end{aligned}
\end{equation}
where $a_1, a_2, a_3, a_4, u$ and $v $ are arbitrary parameters. Further trivial solutions of \eqref{eq5} are obtained by equating both sides of \eqref{eq5} to 0. Solutions that are not trivial will be considered as nontrivial solutions. 

We now  observe that if 
\begin{equation}
(x_1,x_2,x_3,x_4,y_1,y_2,y_3,y_4)=(\alpha_1,\alpha_2,\alpha_3,\alpha_4,\beta_1,\beta_2, \beta_3,\beta_4) \label{solgen}
\end{equation}
 is any solution of the diophantine Eq. \eqref{eq5} and $k_1, k_2$ are any two arbitrary nonzero integers, then 
\begin{multline}
\quad \quad (x_1,x_2,x_3,x_4,y_1,y_2,y_3,y_4)\\
=(k_1\alpha_1,k_1\alpha_2,  k_2\alpha_3,k_2\alpha_4,k_2\beta_1,k_2\beta_2,k_1\beta_3,k_1\beta_4)\quad \quad\quad \quad
\label{solgenk}
\end{multline}
 is also a solution of \eqref{eq5}. It follows that  integer solutions of Eq.~\eqref{eq5} may be obtained from any rational solution of \eqref{eq5} by suitably choosing the integers $k_1$ and $k_2$.

We will now show that the diophantine Eq. \eqref{eq5} is equivalent to a diophantine system consisting of the following three simultaneous  diophantine  equations:
\begin{align}
X_1^5+X_2^5+X_3^5+X_4^5&=Y_1^5+Y_2^5+Y_3^5+Y_4^5, \label{eq5XY}\\
X_1X_2&=Y_1Y_2, \label{eq1X}\\
X_3X_4&=Y_3Y_4. \label{eq1Y}
\end{align}

If $x_i, y_i, i=1,\ldots,\,4$, are any rational numbers satisfying  Eq. \eqref{eq5}, it is readily seen that the rational numbers,
\begin{equation}
\begin{aligned}
X_1&=x_1 x_3, &X_2&=x_2 x_4, &X_3&=-y_1 y_3, &X_4&=-y_2 y_4,\\
Y_1&=-x_1 x_4, &Y_2&=-x_2 x_3, &Y_3&=y_1 y_4, &Y_4&=y_2 y_3,
\end{aligned}
\label{relxyXY}
\end{equation}
satisfy the three equations \eqref{eq5XY}, \eqref{eq1X} and \eqref{eq1Y}. Conversely, if $X_i, Y_i, i=1,\ldots,\,4$, are any rational numbers satisfying the simultaneous equations \eqref{eq5XY}, \eqref{eq1X} and \eqref{eq1Y}, by solving Eqs. \eqref{relxyXY}, we readily obtain rational numbers  $x_i, y_i, i=1,\ldots,\,4$, that satisfy  Eq. \eqref{eq5}. This establishes that the diophantine Eq. \eqref{eq5} is equivalent to the  diophantine system defined by  equations \eqref{eq5XY}, \eqref{eq1X} and \eqref{eq1Y}.

Similarly, it is readily established that the diophantine system defined by  Eq. \eqref{eq5} and the equation,
\begin{equation}
(x_1+x_2)(x_3+x_4)=(y_1+y_2)(y_3+y_4), \label{eq1xy}
 \end{equation}
is equivalent to the diophantine system defined by Eqs. \eqref{eq5XY}, \eqref{eq1X} and \eqref{eq1Y} and the following equation:
\begin{equation}
X_1+X_2+X_3+X_4=Y_1+Y_2+Y_3+Y_4. \label{eq1XY}
\end{equation}

\subsection{A parametric solution}\label{parmsol}
We will now obtain a parametric solution of Eq. \eqref{eq5} by solving the equivalent diophantine system consisting of the simultaneous equations \eqref{eq5XY}, \eqref{eq1X} and \eqref{eq1Y}.

To solve Eqs. \eqref{eq5XY}, \eqref{eq1X} and \eqref{eq1Y}, we write,
\begin{equation}
\begin{aligned}
X_1+X_2&=s_1,\;\;&X_1X_2&=s_2, \;\;&X_3+X_4&=t_1,\;\; &X_3X_4&=t_2,\\
Y_1+Y_2&=S_1,\;\;&Y_1Y_2&=S_2, \;\;&Y_3+Y_4&=T_1,\;\; &Y_3Y_4&=T_2,
\end{aligned}
\label{relXY}
\end{equation}
when Eq. \eqref{eq5XY} may be written as
\begin{multline}
s_1^5 + t_1^5 - 5s_1^3s_2 - 5t_1^3t_2 + 5s_1s_2^2 + 5t_1t_2^2 \\
= S_1^5 + T_1^5 - 5S_1^3S_2 - 5T_1^3T_2 + 5S_1S_2^2 + 5T_1T_2^2, \label{eq5st}
\end{multline}
while equations  \eqref{eq1X} and \eqref{eq1Y} simply reduce to
\begin{equation}
S_2=s_2, \quad T_2=t_2. \label{eq1st}
\end{equation}

We will solve equations  \eqref{eq5st} and \eqref{eq1st} by imposing the following auxiliary condition:
\begin{equation}
T_1=s_1+t_1-S_1. \label{aux1}
\end{equation}

On eliminating $S_2, T_1$ and $T_2$ from Eqs. \eqref{eq5st}, \eqref{eq1st} and \eqref{aux1} we get,
\begin{multline}
s_2^2 - t_2^2 - (s_1^2 + s_1S_1 + S_1^2)s_2 + (s_1^2 - 2s_1S_1 + S_1^2+ 3s_1t_1 - 3t_1S_1 + 3t_1^2)t_2\\
 + (s_1 + t_1)(S_1 - t_1)(s_1^2 - s_1S_1 + S_1^2 + s_1t_1 - t_1S_1 + t_1^2) \label{res1}
\end{multline}

On writing $s_2=t_2+h$, we readily obtain the following solution of Eq. \eqref{res1}:
\begin{equation}
 \begin{aligned}
s_2&=\left[h^2 + \{s_1^2 + (3t_1 - 2S_1)s_1+ 3t_1^2 - 3t_1S_1 + S_1^2 \}h\right.\\
& \quad \quad \left. + (s_1 + t_1)(t_1 - S_1)\{s_1^2+ (t_1 - S_1)s_1 + t_1^2 - t_1S_1 + S_1^2\}\right]\\
& \quad \quad \times \{2h + 3(s_1 + t_1)(t_1 - S_1)\}^{-1},\\
t_2&=\left[-h^2 + (s_1^2 + s_1S_1 + S_1^2)h + (s_1 + t_1)(t_1 - S_1)\{s_1^2 + (t_1 - S_1)s_1 \right.\\
& \quad \quad \left.+ t_1^2 - t_1S_1 + S_1^2\}\right]\{2h + 3(s_1 + t_1)(t_1 - S_1)\}^{-1},
\end{aligned}
\label{vals2t2}
\end{equation} 
where $h$ is an arbitrary parameter.

Since $X_i, Y_i, i=1,\,\ldots, \,4$, satisfy the relations \eqref{relXY}, we will get rational values of $X_i, Y_i$ if and only if $s_1^2-4s_2$, $t_1^2-4t_2$, $S_1^2-4S_2$ and  $T_1^2-4T_2$ are all perfect squares, and we will now choose the parameters suitably to satisfy these four conditions. 

We now write
\begin{equation}
t_1= (S_1^2 - h)/S_1, \label{valt1}
\end{equation}
 and using the values of $s_2, t_2$ given by \eqref{vals2t2}, we get 
\begin{align}
s_1^2-4s_2&=(S_1^2 - S_1s_1 + h)(S_1s_1 - 2h)^2 \nonumber \\
& \quad \quad \times\{S_1^2(S_1^2 + 3S_1s_1 - 3h)\}^{-1} \label{conds12}\\
t_1^2-4t_2&=(S_1^2 - S_1s_1 + h)(S_1^2 + 2S_1s_1 - h)^2  \nonumber \\
& \quad \quad \times\{S_1^2(S_1^2 + 3S_1s_1 - 3h)\}^{-1}. \label{condt12}
\end{align}

It  readily follows from \eqref{conds12} and \eqref{condt12} that both $s_1^2-4s_2$ and $t_1^2-4t_2$ will  become  perfect squares if we choose $s_1$ such that $(S_1^2 - S_1s_1 + h)(S_1^2 + 3S_1s_1 - 3h)$ is a perfect square, and accordingly we choose
\begin{equation}
s_1= \{(3m^2 + 1)h - (m^2 - 1)S_1^2\}/\{(3m^2 + 1)S_1\}, \label{vals1}
\end{equation}
where $m$ is an arbitrary rational parameter. 

Next we will choose our parameters such that $S_1^2-4S_2$ becomes a perfect square.  With the values of $T_1,  s_2, t_2, t_1$  and $s_1$ defined by \eqref{aux1}, \eqref{vals2t2},  \eqref{valt1} and \eqref{vals1} respectively, we get,
\begin{multline}
S_1^2-4S_2=\{(m^2 - 1)(3m^2 + 1)^2h^2 + 2S_1^2(m^4 - 1)(3m^2 + 1)h\\
 + S_1^4m^2(m^2 + 3)^2\}/\{(3m^2 + 1)S_1\}^{2}. \label{condS12}
\end{multline}

Thus, $S_1^2-4S_2$ will become a perfect square if we choose $h$ such that the numerator on the right-hand side of \eqref{condS12} becomes a perfect square. This numerator is a quadratic function of $h$ and we readily find that it becomes a perfect square when we choose
\begin{equation}
\begin{aligned}
h=-2S_1^2u\{(m + 1)(m^2 + 1)u - m(m^2 + 3)\}\\
\quad \quad \times [(3m^2 + 1)\{(m + 1)u^2 - m + 1\}]^{-1}, \label{valh}
\end{aligned}
\end{equation}
where $u$ is an arbitrary rational parameter.

 With the values of $T_1, s_2, t_2, t_1, s_1$  and $h$ defined by \eqref{aux1}, \eqref{vals2t2}, \eqref{valt1}, \eqref{vals1} and \eqref{valh} respectively, we will finally choose a suitable value of $u$ such that $T_1^2-4T_2$ also becomes a perfect square.

We observe that 
\begin{equation}
T_1^2-4T_2=S_1^2\phi(u)/[(3m^2 + 1)\{(m + 1)u^2 - m + 1\}]^2, \label{condT12}
\end{equation}
where
\begin{multline}
\phi(u)=(m^6 - 26m^4 - 31m^2 - 8)(m + 1)^2u^4 - 4m(m + 1)(m^2 + 3)\\
\times (m^4 - 6m^2 - 3)u^3 + 2(m - 1)(m + 1)(3m^6 + 28m^4 + 31m^2 + 2)u^2\\
 - 4m(m - 1)(m^2 + 3)(m^4 + 6m^2 + 1)u + m^2(m + 1)^2(m - 1)^4. \label{valphi}
\end{multline}

It follows from \eqref{condT12} that $T_1^2-4T_2$ will become a perfect square if $u$ is so chosen that $\phi(u)$ is a perfect square. Now $\phi(u)$ is a quartic function of $u$  and,  following a method described by Fermat (as quoted by Dickson \cite[p. 639]{Di}), we obtain    the following value  of $u$ such that $\phi(u)$ becomes a perfect square:
\begin{multline}
u=m(m + 1)^2(m - 1)^3(m^2 + 3)(7m^6 + 23m^4 + 29m^2 + 5)(2m^{14} - 41m^{12} \\
- 328m^{10} - 967m^8 - 1382m^6 - 1047m^4 - 308m^2 - 25)^{-1}. \label{valu}
\end{multline}

We have now chosen the parameters suitably so that $s_1^2-4s_2$, $t_1^2-4t_2$, $S_1^2-4S_2$ and  $T_1^2-4T_2$ are all perfect squares. Using the relations \eqref{relXY}, we now get rational values of $X_i, Y_i,\; i=1,\ldots,4$. We thus obtain a solution of  the simultaneous diophantine equations  \eqref{eq5XY}, \eqref{eq1X} and \eqref{eq1Y} and, on appropriate scaling, this solution may be written as follows:
\begin{equation}
\begin{aligned}
X_1&=(m-1)(m+1)^2f_1(m)f_2(-m),\\
X_2&=-(m+1)(m-1)^2f_1(-m)f_2(m),\\
X_3&=(m-1)^2f_3(m)f_4(m),\\
X_4&=(m+1)^2f_3(-m)f_4(-m),\\
Y_1&=(m-1)^3f_1(m)f_2(m), \\ 
Y_2&=-(m+1)^3f_1(-m)f_2(-m),\\
Y_3&=-(m-1)(m+1)f_3(m)f_4(-m), \\
 Y_4&=-(m-1)(m+1)f_3(-m)f_4(m),
\end{aligned}
\label{parmsol1}
\end{equation}
where the functions $f_i(m),\;i=1,\ldots,4,$ are defined, in terms of an arbitrary parameter $m$,  as follows:
\begin{equation}
\begin{aligned}
f_1(m)&=5m^{14} + 7m^{13} + 71m^{12} + 30m^{11} + 345m^{10} + 17m^9\\
 & \quad + 907m^8 - 60m^7+1311m^6 - 71m^5\\
 & \quad + 1109m^4 + 62m^3 + 323m^2 + 15m + 25,\\
f_2(m)&=m^{10} + 7m^9 + 29m^8 + 44m^7 + 122m^6 + 98m^5 \\
 & \quad+ 202m^4 + 92m^3 + 133m^2 + 15m + 25,\\
f_3(m)&=5m^{14} + 21m^{13} + 29m^{12} + 202m^{11} + 109m^{10} \\
 & \quad+ 755m^9 + 173m^8+ 1388m^7 + 23m^6 + 1259m^5 \\
 & \quad- 177m^4 + 426m^3 - 137m^2 + 45m - 25,\\
f_4(m)&=m^{11} + 6m^{10} + 8m^9 + 57m^8 + 46m^7 + 184m^6 + 92m^5 \\ 
& \quad+ 294m^4 + 89m^3 + 202m^2 + 20m + 25.
\end{aligned}
\label{deffm}
\end{equation}

We note that since we had imposed the auxiliary condition \eqref{aux1}, the solution given by \eqref{parmsol1} also satisfies Eq. \eqref{eq1XY}.

As we have already noted in Section~\ref{prelim}, every solution of the diophantine system defined by the equations equations  \eqref{eq5XY}, \eqref{eq1X} and \eqref{eq1Y} readily yields a solution of the diophantine Eq. \eqref{eq5}. Accordingly, we obtain the following solution of the diophantine Eq. \eqref{eq5}:
\begin{equation}
\begin{aligned}
x_1 &= (m - 1)f_1(m), \quad & x_2 &= (m + 1)f_1(-m), \\ x_3 &= (m + 1)^2f_2(-m), \quad & x_4 &= -(m - 1)^2f_2(m),\\
y_1 &= (m - 1)f_3(m), \quad & y_2 &= (m + 1)f_3(-m), \\ y_3 &= -(m - 1)f_4(m), \quad & y_4 &= -(m + 1)f_4(-m),
\end{aligned}
\label{parmsol1eq5}
\end{equation}
where the functions  $f_i(m),\; i=1,\ldots, 4$, are defined by \eqref{deffm} in terms of the arbitrary parameter $m$. 

Since the values of  
 $X_i, Y_i,\; i=1,\ldots,4$, given by \eqref{parmsol1} satisfy the  additional condition \eqref{eq1XY}, it follows that the values of  $x_i, y_i,\; i=1,\ldots,4$, given by \eqref{parmsol1eq5} satisfy the  additional condition \eqref{eq1xy}.

We note that since any solution \eqref{solgen} of Eq. \eqref{eq5} yields another solution \eqref{solgenk}, we may choose $k_1, k_2$ suitably to derive  a solution of \eqref{eq5} satisfying the additional condition $x_1+x_2=y_1+y_2$. We thus obtain the following solution of  Eq. \eqref{eq5}:
\begin{equation}
\begin{aligned}
x_1 &= (m - 1)f_1(m)f_5(m), & x_2 &= (m + 1)f_1(-m)f_5(m),\\
 x_3 &= (m + 1)^2f_2(-m)f_6(m), & x_4 &= -(m - 1)^2f_2(m)f_6(m),\\
 y_1 &= (m - 1)f_3(m)f_6(m), & y_2 &= (m + 1)f_3(-m)f_6(m),\\
 y_3 &= -(m - 1)f_4(m)f_5(m), & y_4 &= -(m + 1)f_4(-m)f_5(m),
\end{aligned}
\label{parmsol2eq5}
\end{equation}
where the functions  $f_i(m),\; i=1,\ldots,4,$ are defined as before while 
\begin{equation}
\begin{aligned}
f_5(m)&=5m^8 - 12m^6 - 90m^4 - 124m^2 - 35,\\
f_6(m)&=5m^8 + 44m^6 + 94m^4 + 108m^2 + 5,
\end{aligned}
\end{equation}
and $m$ is an arbitrary parameter. The solution \eqref{parmsol2eq5} satisfies the additional conditions,
\begin{align}
x_1+x_2&=y_1+y_2, \label{condxy12}\\
x_3+x_4&=y_3+y_4, \label{condxy34}
\end{align}
The condition \eqref{condxy12} is satisfied, as already mentioned,  by suitable choice of $k_1, k_2$ while the condition \eqref{condxy34}
is automatically satisfied since the solution \eqref{parmsol1eq5} satisfies the condition \eqref{eq1xy}.

A second solution of  Eq. \eqref{eq5}, derived similarly from the solution \eqref{parmsol1eq5} by first choosing $k_1, k_2$ such that $x_1+x_2=y_3+y_4$ and then renaming the variables $y_3, y_4$ as $y_1, y_2$ respectively, is as follows: 
\begin{equation}
\begin{aligned}
x_1 &= (m - 1)f_1(m)f_7(m), \quad &x_2 &= (m + 1)f_1(-m)f_7(m), \\
x_3 &= (m + 1)^2f_2(-m)f_8(m), \quad &x_4 &= -(m - 1)^2f_2(m)f_8(m),\\
y_1 &= -(m - 1)f_4(m)f_8(m),\quad & y_2 &= -(m + 1)f_4(-m)f_8(m),\\
 y_3 &= (m - 1)f_3(m)f_7(m), \quad &y_4 &= (m + 1)f_3(-m)f_7(m), \\
\end{aligned}
\label{parmsol3eq5}
\end{equation}
where the functions  $f_i(m),\; i=1,\ldots,4,$ are defined as before while 
\begin{equation}
\begin{aligned}
f_7(m)&= -m^2 - 1,\\
f_8(m)&=m^6 + 4m^4 + 9m^2 + 2,\\
\end{aligned}
\end{equation}
and, as before, $m$ is an arbitrary parameter. This solution also satisfies the additional conditions \eqref{condxy12} and \eqref{condxy34}.

As a numerical example, when $m=3$, the solution \eqref{parmsol1eq5} yields the following solution of Eqs. \eqref{eq5} and \eqref{eq1xy}:
\begin{equation}
\begin{aligned}
(x_1,x_2,x_3,x_4)&=(35330, 25801, 2407, -1492),\\
(y_1,y_2,y_3,y_4)&= ( -19814,  32807,  1672,  2633).
\end{aligned}
\end{equation}

As a second example, when $m=3$, the solution \eqref{parmsol3eq5} yields the following solution of the simultaneous equations \eqref{eq5}, \eqref{condxy12} and \eqref{condxy34}:
\begin{equation}
\begin{aligned}
(x_1,x_2,x_3,x_4)&=(129005,  176650,  105932, -170897),\\
 (y_1,y_2,y_3,y_4)&= (186943, 118712,  -164035,  99070).
\end{aligned}
\end{equation}

\subsection{More parametric solutions}\label{moreparmsol}

We will now describe a method of generating infinitely many parametric solutions of the diophantine Eq. \eqref{eq5}.

We will follow the method of solution described in Section~\ref{parmsol} till we reach the stage where we need to choose the parameter $u$ such that the quartic function $\phi(u)$, defined by \eqref{valphi}, becomes a perfect square. We will now show that there exist infinitely many values of $u$, given by rational functions of $m$, that will make $\phi(u)$ a perfect square, and accordingly we can obtain infinitely many parametric solutions of Eq. \eqref{eq5}.

The requirement that $\phi(u)$ becomes a perfect square amounts to finding rational points on the curve,
\begin{equation}
v^2=\phi(u). \label{quarticec}
\end{equation}
Since $\phi(u)$ is defined by \eqref{valphi},  Eq. \eqref{quarticec} represents a quartic model of an elliptic curve over the field $\mathbb{Q}(m)$, and by a birational transformation, the quartic curve reduces to the Weierstrass form of the elliptic curve which is as follows:
\begin{multline}
V^2 = U^3 - 432(m - 1)(m + 1)(325m^{10} + 955m^8 + 1266m^6 + 470m^4\\
 + 57m^2 - 1)U - 3456(5m^4 + 2m^2 + 1)(875m^{10} + 2885m^8\\
 + 3822m^6 + 1450m^4 + 183m^2 + 1)(m - 1)^2(m + 1)^2. \label{cubicec}
\end{multline}

The aforementioned birational transformation is defined by
\begin{equation}
\begin{aligned} 
u&=6m(m - 1)\{(m + 1)^2(m - 1)^2U - 420m^{10} - 4812m^8 \\
& \quad \quad - 12648m^6 - 14232m^4 - 4404m^2 - 348\}/\psi(U,V,m),\\
v&=m(m - 1)\{2(m + 1)^3(m - 1)^3U^3 - 36(m - 1)(m + 1)\\
& \quad \quad \times(35m^{10} + 401m^8 + 1054m^6 + 1186m^4 + 367m^2 + 29)U^2\\
& \quad \quad - (m^2 - 1)^3V^2- 864(3m^2 + 1)(m^2 + 3)(5m^4 + 10m^2 + 1)\\
& \quad \quad \times(7m^6 + 23m^4 + 29m^2 + 5)V + 1728(m^2 - 1)^2\\
& \quad \quad \times(42875m^{20}+ 497350m^{18} + 2290155m^{16} + 5717736m^{14}\\
& \quad \quad  + 8360982m^{12} + 7151748m^{10} + 3327950m^8 \\
& \quad \quad+ 821800m^6 + 97551m^4 + 3494m^2 - 89)/\psi^2(U,V,m),
\end{aligned}
\label{biratuv}
\end{equation}
where
\begin{multline}
\psi(U,V,m)=6(m^2 + 3)(m^4 + 6m^2 + 1)U + (m^2 - 1)V \\
- 72(m^2 - 1)(m^2 + 3)(35m^8 + 86m^6 + 108m^4 + 26m^2 + 1),
\end{multline}
and by
\begin{equation}
\begin{aligned}
U&=\{6(m - 1)(m + 1)(3m^6 + 28m^4 + 31m^2 + 2)u^2 \\
& \quad \quad- 36m(m - 1)(m^2 + 3)(m^4 + 6m^2 + 1)u \\
 & \quad \quad+ 18m(m + 1)(m - 1)^2v+ 18m^2(m + 1)^2(m - 1)^4\}/u^2,\\
V&=-108m\{m(m^2 + 3)(m^4 - 6m^2 - 3)(m - 1)^2(m + 1)^2u^3\\
  & \quad \quad- (3m^6 + 28m^4 + 31m^2 + 2)(m + 1)^2(m - 1)^3u^2 \\
 & \quad \quad	+ (m - 1)(m^2 + 3)(m^4 + 6m^2 + 1)uv + 3m(m + 1)(m^2 + 3)\\
 & \quad \quad \times(m^4 + 6m^2 + 1)(m - 1)^3u - m(m + 1)^2(m - 1)^4v\\
 & \quad \quad - m^2(m + 1)^3(m - 1)^6\}/u^3.
\end{aligned}
\label{biratUV}
\end{equation}

We already know a point $P$ on the quartic curve \eqref{quarticec} whose abscissa is given by \eqref{valu}. The point $P^{\prime}=(U_0,V_0)$ on the cubic curve \eqref{cubicec}, corresponding to the  point $P$ on the quartic curve \eqref{quarticec}, is given by
\begin{equation}
\begin{aligned}
U_0&=12(75m^{28} + 1010m^{26} + 11944m^{24} + 103096m^{22} + 585657m^{20}\\
& \quad \quad + 2202226m^{18} + 5635746m^{16} + 10027936m^{14} + 12482909m^{12} \\
& \quad \quad + 10709526m^{10} + 6063588m^8 + 2067944m^6 + 398591m^4\\
& \quad \quad  + 39750m^2 + 1650)/(d^2(m),\\
V_0&=216(125m^{42} + 2525m^{40} + 12350m^{38} - 138015m^{36} - 2822345m^{34}\\
 & \quad \quad - 24701264m^{32} - 140086792m^{30} - 573149148m^{28} - 1776227438m^{26} \\
& \quad \quad - 4275792154m^{24} - 8087224924m^{22} - 12040781858m^{20} \\
& \quad \quad - 14031203010m^{18}- 12641030116m^{16} - 8645319848m^{14}\\
 & \quad \quad- 4384538092m^{12} - 1605427583m^{10} - 411694779m^8 - 71091250m^6\\
 & \quad \quad - 7771895m^4 - 478725m^2 - 12500)/(d^3(m)),
\end{aligned}
\end{equation}
where
\begin{equation}
d(m)=(m^2 + 3)(m - 1)(m + 1)(7m^6 + 23m^4 + 29m^2 + 5).
\end{equation}

When $m=2$, the curve \eqref{cubicec} reduces to
\begin{equation}
V^2 = U^3 - 863202096U - 5268270761856, \label{cubicecspl}
\end{equation}
and a rational point on the  curve \eqref{cubicecspl}  corresponding to the point $P^{\prime}=(U_0, V_0)$ on the curve \eqref{cubicec} is 
\[(3346068693496/43020481, 5630105905921711808/282171334879).\]
Since this rational point on the elliptic curve \eqref{cubicecspl} does not have integer coordinates,  it follows from the Nagell-Lutz theorem \cite[p. 56]{Si} on elliptic curves that this is not a point of finite order.  We can thus find infinitely many rational points on the curve \eqref{cubicecspl} using the group law.

Since in the special case $m=2$, the point on the curve \eqref{cubicecspl} corresponding to the point $P^{\prime}$ is not of finite order, it follows that, for any arbitrary rational value of $m$,  the point $P^{\prime}$ on the curve \eqref{cubicec} cannot  be a point of finite order. We can thus generate infinitely many rational points on the curve \eqref{cubicec} using the group law. Each of these rational points will yield a corresponding rational point on the curve \eqref{quarticec} and  we thus  get infinitely many rational values  of $u$, in terms of the parameter $m$, for which $\phi(u)$  becomes a perfect square, and hence we can  obtain infinitely many parametric solutions of the diophantine Eq. \eqref{eq5}.

The point $P^{\prime}$ on the curve \eqref{cubicec} corresponds to the point $P$ on the quartic curve \eqref{quarticec} and naturally yields the solution \eqref{parmsol1} of \eqref{eq5} already obtained in Section~\ref{parmsol}. The point $2P^{\prime}$  is too cumbersome to write and  yields a parametric solution of Eq. \eqref{eq5} in terms of polynomials  of  degree 74. We accordingly do not give this solution explicitly. However, the above analysis, using elliptic curves, shows that there exist infinitely many parametric solutions of the diophantine Eq. \eqref{eq5} and they can be effectively computed.

\section{Some open problems}\label{openprob} 
Any nontrivial solution of the diophantine equation  \eqref{eq5} in which $y_4=0$ immediately yields a solution in integers of the equation
\begin{equation}
 (x_1^5+x_2^5)(x_3^5+x_4^5)=y_1^5+y_2^5. \label{eq5rev}
\end{equation}
The solutions of \eqref{eq5} obtained in Section~\ref{parmsol} did not yield any such solution. We, however, obtained  four solutions of   Eq. \eqref{eq5rev} by computer trials. The sextuples $(x_1, x_2, x_3, x_4, y_1, y_2)$ giving these four solutions are as follows:
\begin{equation}
\begin{aligned}
&(8, -1, 25, 21, 109, 213),  \quad (19, 12, 6, 4, 41, 119), \quad (2, -1, 77, 83, 136, 174), \\
&(67575, 56763, 21624, -2703, 1556222517, 796376781). 
\end{aligned}
\end{equation}
The first and last of these solutions satisfy the additional condition, 
\[(x_1+x_2)(x_3+x_4)=y_1+y_2.\]
It would of interest to determine whether there exist infinitely many nontrivial solutions of the diophantine Eq. \eqref{eq5rev}.

It would also be interesting to consider the diophantine Eq. \eqref{eqn}
when $n  > 5$. While it is not inconceivable that there exist integer solutions of Eq. \eqref{eqn} when $n=6$, it is unlikely that that there are any solutions of the diophantine equation \eqref{eqn} when $n > 6$.

\noindent Ajai Choudhry, 13/4 A Clay Square, Lucknow - 226001, India

\noindent E-mail address: ajaic203@yahoo.com

\medskip 

\noindent Oliver Couto, Apartment 501,  2166 Lakeshore Road, \\
Burlington, Ontario, L7R 2B6, Canada

\noindent Email: samson@celebrating-mathematics.com
\end{document}